\documentclass [10pt,a4paper] {article}
\usepackage{color}
\usepackage{array}
\usepackage{tabularx}
\usepackage{amsmath}
\usepackage{amsfonts}
\usepackage{amssymb}
\usepackage{bbm}
\usepackage{mathrsfs}
\usepackage{latexsym}
\usepackage{epsfig}
\newtheorem{definition}{Definition}[section]

\newtheorem{lemma}[definition]{Lemma}
\newtheorem{theorem}[definition]{Theorem}
\newtheorem{proposition}[definition]{Proposition}
\newtheorem{corollary}[definition]{Corollary}
\newtheorem{example}[definition]{Example}
\newtheorem{remark}[definition]{Remark}

\begin{document}


\newcommand{\proof}{\paragraph{Proof}}
\newcommand{\fin}{$\Box$\\}
\newcommand{\ds}{\displaystyle}
\newcommand{\saut}[1]{\hfill\\[#1]}
\newcommand{\vsp}{\vspace{.15cm}}
\newcommand{\difrac}{\displaystyle \frac}
\newcommand{\dist}{\textrm{dist}}
\newcommand{\mbf}{\textbf}
\newcommand{\levy}{\mathscr{B}}
\newcommand{\sheet}{\mathbbm{B}}
\newcommand{\sifbm}{\mathbf{B}}
\newcommand{\alphar}{\texttt{\large $\boldsymbol{\alpha}$}}
\title{A set-indexed fractional Brownian motion}
\author{
\begin{minipage}{2.3in}
\centering
Erick Herbin  \vspace{5mm} \\ 
\small Dassault Aviation\\ 
78 quai Marcel Dassault\\
92552 Saint-Cloud Cedex\\
France \\
\small erick.herbin@dassault-aviation.fr
\end{minipage}
\hfill
\begin{minipage}{2.3in}
\centering
Ely Merzbach\thanks{Corresponding author. Tel.: +972-2-9941168}  \vspace{5mm} \\ 
\small Dept. of Mathematics\\
Bar Ilan University\\
52900 Ramat-Gan\\
Israel\\
\small merzbach@macs.biu.ac.il
\end{minipage}
 }
\date{June 2004, May 2005(revised version)}

\maketitle
\begin{abstract}
We define and prove the existence of a fractional Brownian motion indexed by a
collection of closed subsets of a measure space. This process is a
generalization of the set-indexed Brownian motion, when the condition of
independance is relaxed.
Relations with the L\'evy fractional Brownian motion and with the fractional
Brownian sheet are studied.
We prove stationarity of the increments and a property of self-similarity with
respect to the action of solid motions.
Moreover, we show that there no "really nice" set indexed fractional Brownian 
motion other than set-indexed Brownian motion.
Finally, behavior of the set-indexed fractional Brownian motion along increasing paths is
analysed.
\end{abstract}

{\sl AMS classification\/}: 62\,G\,05, 60\,G\,15, 60\,G\,17, 60\,G\,18.

{\sl Keywords\/}: fractional Brownian motion, Gaussian processes, stationarity,
self-similarity, set-indexed processes.


\section{Introduction}
Recently, different developments around fractional Brownian motion were set up
and widely used to describe complex or chaotic phenomena in several fields of
sciences. For instance, let us mention theoretical aspects such as stochastic
calculus with respect to fBm (e. g. \cite{nualart}), 
and some more applied aspects such as its use in finance
(\cite{oksendal}, \cite{sottinen}) or in data traffic modeling (\cite{jlvrr}).
Here we define a new very natural set-indexed generalization of fractional
Brownian motion. The set-indexed fractional Brownian motion studied here seems
to be well-adapted to model problems in applied mathematics (see
\cite{beatrice}).

Fractional Brownian motion was defined by B. B. Mandelbrot and J. W. Van Ness,
and extended essentially into two directions. One is generally called 
multifractional Brownian motion, replacing the index parameter of
self-similarity (called also the Hurst parameter) by a real measurable function
(see \cite{benassi}, \cite{peltier}). The other one are multiparameter fractional Brownian motions in which
the set of the indices is a subset of the Euclidean space $\mathbf{R}^N$ (see
\cite{kamont} and \cite{talagrand}, \cite{xiao} for trajectory properties).

At least two types of multiparameter fractional Brownian motions were
introduced. One is called L\'evy fractional Brownian motion because it extends
L\'evy Brownian motion and the other is called fractional Brownian sheet
because it can be seen as an extension of the Brownian sheet. We refer to
\cite{herbin} for a survey on all these processes.
Moreover, in \cite{herbin}, some multiparameter extensions of 
multifractional Brownian motion are well studied.

The frame of set-indexed processes is not only a new step of generalization of
multiparameter processes, but it provides a real tool in modeling (e. g.
\cite{giem}). Then a set-indexed extension of fractional Brownian motion is 
hoped to provide a powerful model for multidimensional phenomena.  
The set-indexed fractional Brownian motion introduced here is a simple extension
of set-indexed Brownian motion, also called white noise (see \cite{adlerBM},
\cite{davar}, \cite{goldie} and \cite{pyke} for studies and applications), 
and possesses the main properties required for fractional Brownian motion. 
Moreover, by choosing the class of sets parametrizing the process, we get great
flexibility and possibilities to obtain particular types of fractional Brownian
motion.

In this paper, we prove that our definition of set-indexed fractional Brownian
motion satisfies both self-similarity and a condition of stationarity. Moreover
such a process is a time-changed classical fractional Brownian motion on each
increasing path (flow).\\
Conversely, we compute the covariance function of any self-similar and
stationary set-indexed process. For Gaussian processes, we get an extension of
our set-indexed fractional Brownian motion.

This paper is divided as follow :\\
In the next section, we present the general framework needed for
set-indexed processes. We define the concept of set-indexed fractional
Brownian motion (sifBm). We prove existence of this process showing that its
covariance function is positive definite. Moreover, we compare our
definition to previous definitions given in some particular cases and our
definition seems to be quite natural and satisfactory. 
The two fractal properties which are
stationarity and self-similarity are studied in section 3. 
As it will be see, stationarity of
increments can be defined in different non equivalent ways.
In section 4, we discuss the possibility of finding a characterization of  
set-indexed fractional Brownian motion by the two main properties: stationarity
and self-similarity. Moreover, we show that there no "really nice" set indexed 
fractional Brownian motion other than set-indexed Brownian motion.
More precisely, theorem \ref{thCstat} states that the only set-indexed Gaussian
process satisfying 
$E\left[\left(\Delta X_C\right)^2\right]=m(C)^{2H}$ ($\forall
C\in\mathcal{C}$), where $m$ is the Lebesgue measure, $\Delta X$ the increment
process, and $H\in(0,1)$, is the white noise ($H=1/2$).
Section 5 deals with the problem of continuity. In fact,
we show that sifBm is continuous when set-indexed Brownian motion
is also continuous. 
Finally in the last section, we study sifBm on increasing paths.

\section{Framework and definition}

\subsection{Indexing collection, set-indexed processes}
Let $\mathcal{T}$  be a  locally compact complete separable metric and measure 
space with metric $d$ and measure $m$.
All processes will be indexed by a class $\cal{A}$ of compact connected subsets
of $\mathcal{T}$.

In what follows, for any class of sets ${\cal D}$, the class of finite unions of
sets from ${\cal D}$ will be denoted by ${\cal D}(u)$.  In the terminology of
\cite{Ivanoff}, we assume that ${\cal A}$ is an {\it indexing collection}:

\begin{definition}\label{basic}
A nonempty class ${\cal A}$ of compact, connected subsets of $\mathcal{T}$ is 
called an {\em indexing collection}\index{indexing collection} 
if it satisfies the following:
\begin{enumerate}
 \item $\emptyset\in {\cal A}$, and $A^{\circ}\neq A$ 
if $A\notin\left\{ \emptyset, \mathcal{T} \right\}$.
In addition, there exists an increasing sequence 
$\left(B_n\right)_{n\in\textbf{N}}$ of sets in $\mathcal{A}(u)$ such that
$\mathcal{T}=\bigcup_{n\in\textbf{N}}B^{\circ}_n$.

\item ${\cal A}$ is closed under arbitrary intersections and if 
$A,B\in {\cal A}$ are nonempty, then $A\cap B$ is nonempty. 
If $(A_i)$ is an increasing sequence in ${\cal A}$ and if there exists
$n\in\textbf{N}$ s. t. for all $i$, $A_i\subseteq B_n$ 
then $\overline{\bigcup_i A_i}\in {\cal A}$.

\item The $\sigma$-algebra generated by ${\cal A}$, $\sigma ({\cal A})={\cal B}$,
the collection of all Borel sets of $\mathcal{T}$.

\item {\em Separability from above}\\
There exists an increasing sequence of finite subclasses 
${\cal A} _n=\{A_1^n,...,A_{k_n}^n\}$ of ${\cal A}$ closed under intersections 
and satisfying $\emptyset, B_n \in\mathcal{A}_n(u)$  
and a sequence of functions 
$g_n:\mathcal{A}\rightarrow\mathcal{A}_n(u)\cup\left\{\mathcal{T}\right\}$ 
such that
\begin{enumerate}
\item $g_n$ preserves
arbitrary intersections and finite unions (i.e.\linebreak 
$g_n(\bigcap_{A\in {\cal A}'}A )=\bigcap_{A\in {\cal A}'}g_n(A )$ 
for any ${\cal A}'\subseteq {\cal A}$, and if 
$\bigcup_{i=1}^kA_i=\bigcup_{j=1}^mA_j'$, then
$\bigcup_{i=1}^kg_n(A_i)=\bigcup_{j=1}^mg_n(A_j')$),
\item
 for each $A\in {\cal A}$, $A\subseteq (g_n(A))^{\circ}$,
 \item
$g_n(A)\subseteq g_m(A)$ if $n\geq m$,
\item  for each $A\in{\cal A}$, $A=\bigcap_ng_n(A)~ $,
\item if $A,A'\in{\cal A}$ then for every $n$, $g_n(A) \cap A' \in {\cal A}
$, and if
$A'\in{\cal A} _n$ then $g_n(A) \cap A' \in {\cal A} _n$.
\item   $g_n(\emptyset)=\emptyset~\forall n$.
\end{enumerate}

\item Every countable intersection of sets in ${\cal A} (u)$ may be expressed as
the closure of a countable union of sets in ${\cal A}$.

\end{enumerate}
 (Note: ` $\subset$' indicates strict
inclusion and `$\overline{(\cdot)}$' and`$(\cdot )^{\circ}$' denote
respectively the closure and the interior of a set.)

\end{definition}

We shall define the semi-algebra\index{semi-algebra} ${\cal C}$
\index{${\cal C}$} to be the class of all subsets of $\mathcal{T}$ of the form
$$C=A\setminus B,~A\in {\cal A},~ B\in {\cal A} (u).$$
$\mathcal{C}$ is closed under intersections and any set in $\mathcal{C}(u)$
\index{${\cal C}(u)$}  
may be expressed as a finite  disjoint union of sets in $\mathcal{C}$. 
Note that if
$B=\bigcup_{i=1}^kA_i \in {\cal A} (u)$, without loss of generality we can 
require that for each $i$, $A_i\not\subseteq \bigcup_{j\ne i}A_j$. 
Such a representation of $B\in {\cal A} (u)$ will be called
{\em extremal}\index{extremal representation}. 
If $C=A\setminus B,~A\in {\cal A},~ B\in {\cal A} (u)$, 
then the representation of $C$ is called extremal if that of $B$ is. 
Unless otherwise stated, it will always be assumed that all
representations of sets in ${\cal A} (u)$ and $C$ are extremal. 

Numerous examples of topological spaces $\mathcal{T}$ and indexing collections ${\cal A}$
satisfying the preceding assumptions may be found in \cite{Ivanoff}. 
In particular,
our framework generalizes the usual multiparameter setting: 
if $\mathcal{T}=\mathbf{R}^N_+$,
then the class 
${\cal A} = \left\{[0,t]:t\in \mathbf{R}^N_+\right\} \cup
\left\{\emptyset\right\}$ 
satisfies all the assumptions.
More generally, we can allow ${\cal A}$ to consist of all the {\em lower layers}
of $\mathbf{R}^N_+$: 
a set $A$ is a lower layer if $[0,t]\subseteq A,~\forall t\in A$.

Now, let $(\Omega,\mathcal{F},P)$ be any complete probability space. 
A filtration (indexed by ${\cal A}$) is a class 
$ \{ { \cal F } _A :A\in{\cal A} \}$  of complete sub-$\sigma$-fields of
$\mathcal{F}$
\index{filtration}\index{${\cal F}_A$}  
which satisfies the following conditions:
\begin{itemize}

\item $ \forall A , B \in { \cal A } $, $ { \cal F } _ A \subseteq {
\cal F } _ B $ , if $ A \subseteq B $.

\item Monotone outer-continuity:
\index{monotone outer-continuous!filtration}
$ {\cal F}_{\bigcap A_i } = \bigcap {\cal F}_{A_i}$ for any decreasing sequence
$(A_i)$ in $ { \cal A }$.

\end{itemize}

\begin{definition} A (${\cal A}$ -indexed) stochastic process\index{stochastic!process,
{\em see process}}\index{process!{\cal A} -indexed}
$X=\{X_A:A\in{\cal A}\}$  is a collection of random variables indexed by
${\cal A}$ with $X_{\emptyset}=0$, and is said to be
adapted\index{process!adapted}\index{adapted process!${\cal A}$-indexed} if
$X_A$ is $\mathcal{F}_A$-measurable,  for every $A\in {\cal A}$. 
$X$ is said to be integrable
\index{process!integrable}\index{integrable process}
(square integrable)
\index{process!square integrable}\index{square integrable process}
if $E\left[|X_A|\right]<\infty$
 ($E\left[(X_A)^2\right]<\infty$) for every $A\in{\cal A}$.
\end{definition}

\begin{remark}\label{rkmultset}
Any multiparameter process $\tilde{X}$ can be considered as a set-indexed 
process $X$, setting $\mathcal{T}=\mathbf{R}^N_{+}$,
$\mathcal{A}=\left\{[0,t];\;t\in\mathbf{R}^N_{+}\right\}$ and
$\tilde{X}_t=X_{[0,t]}$.
\end{remark}

\begin{definition}\label{additive process}
\index{process!additive}\index{additive!process} 
A ($\mathcal{A}$-indexed) stochastic process $X$ is {\em additive} 
if it has an (almost sure) additive extension to $\mathcal{C}(u)$: 
i.e., $X_{\emptyset}= 0$ and if
 $C,C_1, C_2 \in \mathcal{C} (u) $ with $C=C_1\cup C_2$ and 
 $C_1\cap C_2 = \emptyset$, then almost surely
 $$X_{C_1}+X_{C_2}=X_C.$$
\end{definition}

In addition to assuming that  $X_{\emptyset}= 0$, to avoid technicalities we
will generally assume as well that $X_{\emptyset '}= 0$, 
where $\emptyset ':=\bigcap_{A\in{\cal A}}A$.

\begin{definition}\label{increasing} An additive process $X$ is increasing if $X_C\geq 0
~\forall C\in\mathcal{C}$ and if for any decreasing sequence $(A_n)$ in
$\mathcal{A}$,
$X(\bigcap_n A_n)=\lim_n X(A_n)$.
\end{definition}

It is observed in Corollary 1.4.11 of \cite{Ivanoff} that an increasing process in
fact defines a measure on $\mathcal{B}$ for each $\omega \in \Omega$.

\begin{definition}
Let $\Lambda$ be a non-negative increasing function defined on $ {\cal A}$ 
with $\Lambda_{\emptyset}=0$.
We say that an ${\cal A}$-indexed additive process $X$ is a Brownian
motion\index{Brownian motion} with variance measure
\index{Brownian motion!variance measure}
$\Lambda$ if $X_{\emptyset}= 0$, and if  for disjoint sets
$C_1,\dots,C_n\in\mathcal{C}, X_{C_1},\dots,X_{C_n}$ are independent mean-zero 
Gaussian random variables with variances $\Lambda_{C_1},...\Lambda_{C_n}$, 
respectively. 
\end{definition}

\subsection{Increments of a set-indexed process}
The notion of increments for a set-indexed process 
$X=\left\{X_U ;\; U\in\mathcal{A} \right\}$ is not as simple as in the case of
real indices, where it is only the difference between values of the process.

In the case of multiparameter processes, we use to define the increment between
$s,t\in\mathbf{R}^N_{+}$ by
\begin{equation}\label{multinc}
\Delta X_{s,t}=\sum_{r\in\left\{0,1\right\}^N} (-1)^{N-\sum_i r_i}
X_{[s_i+r_i(t_i-s_i)]_i}
\end{equation}
which is different from the simple difference $X_t-X_s$ (see \cite{herbin}).

In the case of set-indexed processes, the increments are defined from
the collection of subsets $\mathcal{C}$.\\
For all $C=U\setminus\bigcup_{1\leq i\leq n} U_i$, we define the increment of
the process $X$ on $C$ by
\begin{equation}
\Delta X_C=X_U-\sum_{i=1}^n X_{U\bigcap U_i}
+\sum_{i<j}X_{U\bigcap (U_i\bigcap U_j)} -\dots
+(-1)^n X_{U\bigcap \left(\bigcap_{1\leq i\leq n}U_i\right)}
\end{equation}
According to remark \ref{rkmultset}, this expression, applied to the multiparameter case, gives the definition 
(\ref{multinc}) of the increments.\\
In the following, it would be important to consider the particular increments
$\mathcal{C}_0=\left\{ C=U\setminus V;\;U,V\in\mathcal{A} \right\}$.
Moreover the definition of the increment process $\Delta X$ can be extended to 
$\mathcal{C}(u)$, the finite unions of elements of $\mathcal{C}$. 
Particularly, for all $U,V\in\mathcal{A}$, 
$\Delta X_{U\bigtriangleup V}$ is well-defined.

\begin{remark}
The process $\Delta X$ could be seen as an 
extension of the process $X$ for the set of indices $\mathcal{C}$.
For all $A\in\mathcal{A}\subset\mathcal{C}$, we have
$\Delta X_A=X_A$ (because $A=A\setminus\emptyset$).\\
In the case of an additive process, the definition of the increment
$\Delta X_C$ coincides with the additive extension of $X_C$ to 
$C\in\mathcal{C}$.
However, if $X$ is not additive, which is the case of a direct definition of the
process for the set of indices $\mathcal{C}$, in general
$\Delta X_C\neq X_C$ for $C\in\mathcal{C}$.
For this reason, we use a different notation for the increments of $X$.
\end{remark}

\begin{remark}\label{RemIncGauss}
If $X=\left\{ X_U;\;U\in\mathcal{A} \right\}$ is Gaussian, then
$\Delta X=\left\{ \Delta X_C;\;C\in\mathcal{C} \right\}$ is clearly 
Gaussian.
\end{remark}

\subsection{Definition of sifBm}

Recall that the fractional Brownian motion $B^H$ is defined to be a mean-zero 
Gaussian process such that
\[
\forall s,t\in\mathbf{R}_{+};\quad 
E\left[\left(B^H_t-B^H_s\right)^2\right]=|t-s|^{2H}
\]
The natural set-indexed extension of this process is to substitute the term 
$|t-s|^{2H}$ with $d(U,V)^{2H}$, where $d$ is some distance between 
two subsets of $\mathcal{T}$.
In this paper, we consider the choice of 
$d(U,V)=m(U\bigtriangleup V)$, where $\bigtriangleup$ is the
symmetric difference between two sets and $m$ is a measure on 
$\mathcal{T}$.

\begin{lemma}\label{deflemm}
Let $m$ be a finite measure on $\mathcal{T}$.
For all $\alpha\in (0,1]$, the function
\[
(U,V)\mapsto m(U)^{\alpha}+m(V)^{\alpha}-m(U\bigtriangleup V)^{\alpha}
\]
is positive definite.
\end{lemma}

\proof
For all measurable subset $U$ of $\mathcal{T}$ such that $m(U)<+\infty$, 
we define the elementary function $f=\mathbbm{1}_U$. Finite linear combinations
of elementary functions are called simple functions.
It is well-known that simple functions are dense in $L^2(m)$.\\
Moreover, for all $U, V \subseteq\mathcal{T}$, we have
\begin{align*}
\mathbbm{1}_{U\bigtriangleup V}&=\mathbbm{1}_{(U\setminus V) \cup (V\setminus U)}\\
&=\mathbbm{1}_U (1-\mathbbm{1}_V) + \mathbbm{1}_V (1-\mathbbm{1}_U)\\
&=\mathbbm{1}_U + \mathbbm{1}_V -2\; \mathbbm{1}_U.\mathbbm{1}_V\\
&=\left(\mathbbm{1}_U - \mathbbm{1}_V\right)^2\\
&=\left|\mathbbm{1}_U - \mathbbm{1}_V\right|
\end{align*}

Then we only have to show that the function
\begin{eqnarray*}
L^2(m)\times L^2(m)&\rightarrow& \mathbf{R}\\
(f,g)&\mapsto& m(f^2)^{\alpha}+m(g^2)^{\alpha}-m(|f-g|^2)^{\alpha}
\end{eqnarray*}
is positive definite.

Let $f_1, f_2, \dots, f_n \in L^2(m)$ and $u_1, u_2, \dots, u_n \in\mathbf{R}$.
We have to show that
\begin{equation}\label{defpos1}
\sum_{i=1}^n\sum_{j=1}^n \left\{ m(f_i^2)^{\alpha}+m(f_j^2)^{\alpha}
-m(|f_i-f_j|^2)^{\alpha} \right\} u_i u_j \geq 0
\end{equation}
Setting $u_0=-\sum_{i=1}^n u_i$ and $f_0=\mathbbm{1}_{\emptyset}$, we can write
\begin{equation}
\sum_{i=1}^n\sum_{j=1}^n \left\{ m(f_i^2)^{\alpha}+m(f_j^2)^{\alpha}
-m(|f_i-f_j|^2)^{\alpha} \right\} u_i u_j 
=-\sum_{i=0}^n\sum_{j=0}^n m(|f_i-f_j|^2)^{\alpha} u_i u_j
\end{equation}
But, for all $\lambda>0$, we have
\begin{align}
\sum_{i=0}^n\sum_{j=0}^n e^{-\lambda m(|f_i-f_j|^2)^{\alpha}} u_i u_j
&=\sum_{i=0}^n\sum_{j=0}^n \left(e^{-\lambda m(|f_i-f_j|^2)^{\alpha}}-1\right) u_i
u_j \nonumber\\
&=-\lambda \sum_{i=0}^n\sum_{j=0}^n m(|f_i-f_j|^2)^{\alpha} u_i u_j
+o(\lambda)
\end{align}
Then, (\ref{defpos1}) is equivalent to 
\begin{equation}\label{defpos2}
\sum_{i=0}^n\sum_{j=0}^n e^{-\lambda m(|f_i-f_j|^2)^{\alpha}} u_i u_j \geq 0
\end{equation}
in the neighborhood of $\lambda=0$.

In $L^2(m)$, let us define the bilinear form $<f,g>=m(fg)$.
If we identify the elements of $L^2(m)$ that are almost everywhere equal, then
$L^2(m)$ is a complete separable metric space with this scalar product.\\
Let us show that there exists a random variable $X$ taking its values in 
$L^2(m)$, such that
\begin{equation}\label{caract}
\forall f\in L^2(m);\quad E\left[e^{i <f,X>}\right]=e^{-\lambda\|f\|^{2\alpha}}
\end{equation}
Consider a $\alpha$-stable real random variable
\[
Y \stackrel{(d)}{=} S_{\alpha}
\left(\left(\cos\frac{\pi\alpha}{2}\right)^{\frac{1}{\alpha}},1,0\right)
\]
where $0<\alpha<1$ 
(see \cite{taqqu}, prop. 1.2.12, p 15).\\
As the function $f\mapsto e^{-\frac{1}{2}\|f\|^2}$ is positive definite
(see remark \ref{remisom}),
by a theorem of Bochner-Minlos (see \cite{hkps} and \cite{kuo}), 
there exists a random variable $G$, such that
\[
\forall f\in L^2(m);\quad E\left[e^{i <f,G>}\right]=e^{-\frac{1}{2}\|f\|^2}.
\]
Moreover, we can suppose that $G$ is independant from $Y$.\\
For all $f\in L^2(m)$, we compute
\begin{align*}
E\left[e^{i <f,Y^{1/2}G>}\right]
&=E\left\{E\left[e^{i <f,Y^{1/2}G>}\mid Y\right]\right\}\\
&=E\left\{E\left[e^{i Y^{1/2} <f,G>}\mid Y\right]\right\}\\
&=E\left[e^{-\frac{Y}{2} \|f\|^2}\right]
\end{align*}
As $E\left[e^{-\gamma Y}\right]=e^{-\gamma^{\alpha}}$ for all $\gamma>0$, 
we get
\begin{equation*}
E\left[e^{i <f,Y^{1/2}G>}\right]
=e^{-2^{-\alpha} \|f\|^{2\alpha}}
\end{equation*}
Then 
$X=\left\{\begin{array}{ll}
\sqrt{2\lambda^{1/\alpha}}\; Y^{1/2}G&\textrm{if }\alpha\neq 1\\
\sqrt{2\lambda}\; G&\textrm{if }\alpha=1
\end{array}\right.$ 
satisfies (\ref{caract}), and 
$f\mapsto e^{-\lambda\|f\|^{2\alpha}}$ is non-negative definite.
That proves (\ref{defpos2}) and the result follows.
\fin

\begin{remark}\label{remisom}
For all family $(f_1,\dots,f_n)$ of $L^2(m)$ provided with the scalar product
previously defined, there exist $p\leq n$ and a family $(x_1,\dots,x_n)$ of
$\mathbf{R}^p$ such that
\begin{equation}\label{isom}
\forall i,j;\quad \|x_i - x_j\|_{\mathbf{R}^p} = \|f_i - f_j\|_{L^2(m)}
\end{equation}
\end{remark}

To show this result, let us consider an orthonormal basis $(e_1,\dots,e_p)$ of 
$Vect(f_i)_{1\leq i\leq n}$, and the canonical basis 
$(\epsilon_1,\dots,\epsilon_p)$ of $\mathbf{R}^p$.
For all $i$, there exists a family $(\lambda^i_1,\dots, \lambda^i_p)$ s.t.
$f_i=\sum_k \lambda^i_k.e_k$. Then the vectors $(x_1,\dots,x_n)$ defined by
$x_i=\sum_k \lambda^i_k.\epsilon_k$ satisfy (\ref{isom}).
\fin
This remark allows to show directly (\ref{defpos1}) using lemma 2.10.8 in
\cite{taqqu}, for $2\alpha\in (0,2]$ i. e. $\alpha\in (0,1]$.

\vsp
Since the existence of a mean-zero Gaussian process is equivalent to the
positive definite property of its covariance function, we can define

\begin{definition}\label{defsifBm}
A mean-zero Gaussian process $\sifbm^H=\left\{\sifbm^H_U ;\; U\in\mathcal{A} \right\}$
satisfying
\begin{equation}\label{defcov}
E\left[{\sifbm^H_U} {\sifbm^H_V}\right]=\frac{1}{2}\left[m(U)^{2H}+m(V)^{2H}
-m(U\bigtriangleup V)^{2H}\right]
\end{equation}
where $H\in (0,1/2]$,
is called a set-indexed fractional Brownian motion (sifBm).
$H$ is the index of self-similarity of the process. 
\end{definition}

Lemma \ref{deflemm} only shows that the right side of equation (\ref{defcov}) is positive definite, which restricts the definition of the sifBm for $H\in (0,1/2]$. Even in the simple case of $\mathcal{A}=\left\{[0,t];\;t\in\mathbf{R}^2_+\right\}\cup\left\{\emptyset\right\}$, some examples can be found where $H>1/2$ leads to a non positive definite expression (\ref{defcov}).

\begin{remark}
If $H=\frac{1}{2}$, the process $\sifbm^{\frac{1}{2}}$ is the well known set-indexed
Brownian motion.
Indeed, let us compute the covariance function of this process
\[
E\left[{\sifbm^{\frac{1}{2}}_U} {\sifbm^{\frac{1}{2}}_V}\right]=\frac{1}{2}
\left[m(U)+m(V)-m(U\bigtriangleup V)\right]
\]
As 
\begin{align*}
m(U\bigtriangleup V)&=m(U\setminus V)+m(V\setminus U)\\
&=m(U\setminus U\cap V)+m(V\setminus U\cap V)\\
&=m(U)+m(V)-2\;m(U\cap V)
\end{align*}
we have
\[
E\left[{\sifbm^{\frac{1}{2}}_U} {\sifbm^{\frac{1}{2}}_V}\right]=m(U\cap V)
\]
which is the covariance function of the set-indexed Brownian motion.
\end{remark}

\begin{remark}
In the case of $\mathcal{T}=\mathbf{R}_{+}$ and $\mathcal{A}=\left\{[0,t];\;
t\in\mathbf{R}_{+}\right\} \cup \left\{\emptyset\right\}$, 
the process $\sifbm^H$ is the classical fBm.
Indeed, the covariance function is
\[
E\left[{\sifbm^H_{[0,s]}} {\sifbm^H_{[0,t]}}\right]=\frac{1}{2} \left[
s^{2H}+t^{2H}-|t-s|^{2H}\right]
\]
which is the covariance function of the fractional Brownian motion.
\end{remark}

\begin{remark}
In the case of $\mathcal{T}=\mathbf{R}_{+}^N$ and $\mathcal{A}$ is the set of
rectangles of the form $[0,t]$, the process $\sifbm^H$ can be seen as a
multiparameter process. Then it is interesting to compare it with the other
known multiparameter fractional Brownian motions (see \cite{herbin}).
\begin{itemize}
\item the L\'evy fractional Brownian motion is a mean-zero Gaussian process 
$\levy^H=\left\{\levy^H_t;\;t\in\mathbf{R}^N_{+}\right\}$ such that 
\begin{equation*}
\forall s,t\in\mathbf{R}^N_{+};\quad
E\left[{\levy^H_s}\; {\levy^H_t}\right]=\frac{1}{2}
\left[\|s\|^{2H}+\|t\|^{2H}-\|t-s\|^{2H}\right]
\end{equation*}
where $H\in(0,1)$.

\item the fractional Brownian sheet is a mean-zero Gaussian process 
$\sheet^H=\left\{\sheet^H_t;\;t\in\mathbf{R}^N_{+}\right\}$ such that 
\begin{equation*}
\forall s,t\in\mathbf{R}^N_{+};\quad
E\left[{\sheet^H_s}\; {\sheet^H_t}\right]=\frac{1}{2} \prod_{i=1}^N
\left[s_i^{2H_i}+t_i^{2H_i}-|t_i-s_i|^{2H_i}\right]
\end{equation*}
where $H=(H_1,\dots,H_N)\in(0,1)^N$.
\end{itemize}
\vsp

As for all $s,t\in\mathbf{R}^N_{+}$,
\begin{equation*}
E\left[{\sifbm^H_{[0,s]}} {\sifbm^H_{[0,t]}}\right]=\frac{1}{2} 
\left[ \prod_{i=1}^N s_i^{2H} + \prod_{i=1}^N t_i^{2H}
-\left( \prod_{i=1}^N s_i + \prod_{i=1}^N t_i - 2\prod_{i=1}^N s_i\wedge t_i \right)^{2H} \right]
\end{equation*}
we can see that, if $m$ is the Lebesgue measure and $N>1$, 
$\sifbm^H$ is different from the two processes $\levy^H$ and $\sheet^H$.
\end{remark}

This fact will be also shown in the next sections in the study of properties of 
sifbm, and its restriction on flows.
It is therefore natural to wonder if the L\'evy fBm and the fBs can have
set-indexed extension. The answer seems to be negative.\\
Actually the definition of the sheet is strongly associated with 
the Euclidean structure of $\mathbf{R}^N$. 
Therefore it is incompatible with a set-indexed viewpoint.\\
Moreover the L\'evy fBm can be seen as a simple one parameter process where the
increment between two points only depends from distance between them.

\section{Fractal properties}
The fractional Brownian motion $B^H$ has two important properties which make it
the most natural fractal process:
\begin{itemize}
\item its increments are stationary
\begin{equation*}
\forall h\in\mathbf{R}_{+};\quad 
\left(B^H_{t+h}-B^H_h\right)_{t\in\mathbf{R}_{+}} \stackrel{(d)}{=} 
\left(B^H_t-B^H_0\right)_{t\in\mathbf{R}_{+}}
\end{equation*}

\item it is self-similar
\begin{equation*}
\forall a\in\mathbf{R}_{+};\quad
\left(B^H_{at}\right)_{t\in\mathbf{R}_{+}} \stackrel{(d)}{=}
\left(a^H B^H_{t}\right)_{t\in\mathbf{R}_{+}}
\end{equation*}

\end{itemize}
Moreover, the fBm is the only Gaussian process which has these two properties.

In this section, we show that in some sense these properties still hold for the set-indexed
fractional Brownian motion. Moreover they characterize the covariance function of
the process between two sets $U$ and $V$ such that $U\subseteq V$.

\subsection{Stationarity of the increments}
Stationarity of increments of a set-indexed process can be defined in
various ways. The set-indexed Brownian motion satisfies all of them, but the
different extensions of {\it fractional} Brownian motion do not.

In the case of $\mathcal{T}=\mathbf{R}_{+}^N$ and 
$\mathcal{A}$ is the collection of rectangles, the classical definition of
stationarity of increments can be studied.\\
A process $X=\left\{ X_{[0,t]};\;t\in\mathbf{R}^N_{+} \right\}$ 
is said to have {\it stationary increments against translations}
if for all $\tau\in\mathbf{R}^N_{+}$, the two processes 
$\left\{ \Delta X_{[\tau,t+\tau]};\;t\in\mathbf{R}^N_{+} \right\}$
and $\left\{ \Delta X_{[0,t]};\;t\in\mathbf{R}^N_{+} \right\}$
have the same law.\\
Both L\'evy fractional Brownian motion and fractional Brownian sheet satisfy this
property of stationarity (see \cite{herbin}).

\begin{remark}
This definition is weaker than stationarity of increments against 
isometries of $\mathbf{R}^N_{+}$, i.e.
for all $g\in\mathcal{G}(\mathbf{R}^N)$, 
\[
\left\{ \Delta X_{g([0,t])};\;t\in\mathbf{R}^N_{+} \right\} \stackrel{(d)}{=}
\left\{ \Delta X_{[0,t]};\;t\in\mathbf{R}^N_{+} \right\}
\] 
where $\mathcal{G}(\mathbf{R}^N)$ is the group of isometries of 
$\mathbf{R}^N$.\\
On the contrary to stationarity against translations, the context of set-indexed
processes imposes additional assumptions to the strict context of multiparameter
processes. 
Actually as the image of $C\in\mathcal{C}$ by any isometry of $\mathbf{R}^N$ 
does not necessarily belong to $\mathcal{C}$, a stability assumption is needed
for the definition to make sense.\\
However in the strict context of multiparameter processes, this 
assumption does not need to be considered. 
The L\'evy fractional Brownian motion satisfy this property of increment
stationarity in the strong sense (see \cite{taqqu}).
\end{remark}

However in general, there is no reason that the sifBm possesses the stationarity
increments property against translations. In some particular cases, we can show
directly using the next lemma that this property is not satisfied.

\begin{lemma}\label{statRNlem}
Let $\sifbm^H=\left\{ \sifbm^H_U;\; U\in\mathcal{A} \right\}$ be a sifBm 
of index $H\in (0,1/2]$.
For all $C=U\setminus \left(\bigcup_{1\leq i\leq n}U_i\right)\in\mathcal{C}$,
where $\forall i\in\left\{1,\dots,n\right\}; U_i\subset U$, we have
\begin{multline}
E\left[\left( \Delta\sifbm^H_C \right)^2\right]=
-\sum_{k=1}^n (-1)^k \sum_{i_1<\dots<i_k} 
m\left(U\bigtriangleup\left[\bigcap_{p\in\left\{i_1,\dots,i_k\right\}}U_p\right]\right)^{2H}\\
+ \frac{1}{2}\sum_{k,l} (-1)^{k+l}
\sum_{i_1<\dots<i_k \atop j_1<\dots<j_l}
m\left(\left[\bigcap_{p\in\left\{i_1,\dots,i_k\right\}}U_p\right]\bigtriangleup
\left[\bigcap_{p\in\left\{j_1,\dots,j_l\right\}}U_p\right]\right)^{2H}
\end{multline}
\end{lemma}

\proof
By definition, 
\begin{equation*}
\Delta\sifbm^H_C=\sifbm^H_U+\sum_{k=1}^n (-1)^k \sum_{i_1<\dots<i_k}
\sifbm^H_{\bigcap_{p\in\left\{i_1,\dots,i_k\right\}}U_p}
\end{equation*}
Then,
\begin{multline*}
E\left[\left(\Delta\sifbm^H_C\right)^2\right]=E\left[\left(\sifbm^H_U\right)^2\right]
+2\sum_k (-1)^k \sum_{i_1<\dots<i_k}
E\left[\sifbm^H_U.\sifbm^H_{\bigcap_{p\in\left\{i_1,\dots,i_k\right\}}U_p}\right]\\
+\sum_{k,l} (-1)^{k+l} \sum_{i_1<\dots<i_k \atop j_1<\dots<j_l}
E\left[\sifbm^H_{\bigcap_{p\in\left\{i_1,\dots,i_k\right\}}U_p}.
\sifbm^H_{\bigcap_{p\in\left\{j_1,\dots,j_l\right\}}U_p}\right]
\end{multline*}
and, using the covariance function of $\sifbm^H$
\begin{equation}\label{varinterm}
\begin{split}
E\left[\left(\Delta\sifbm^H_C\right)^2\right]&=m(U)^{2H}\\
+&\sum_k (-1)^k\sum_{i_1<\dots<i_k}\left\{m(U)^{2H}
+m\left(\bigcap_{p\in\left\{i_1,\dots,i_k\right\}}U_p\right)^{2H}
-m\left(U\bigtriangleup\left[
\bigcap_{p\in\left\{i_1,\dots,i_k\right\}}U_p\right]\right)^{2H}
\right\}\\
+&\frac{1}{2}\sum_{k,l} (-1)^{k+l}
\sum_{i_1<\dots<i_k \atop j_1<\dots<j_l}
\left\{ m\left(\bigcap_{p\in\left\{i_1,\dots,i_k\right\}}U_p\right)^{2H}
+m\left(\bigcap_{p\in\left\{j_1,\dots,j_l\right\}}U_p\right)^{2H} \right.\\
&\phantom{+\frac{1}{2}\sum_{k,l} (-1)^{k+l}
\sum_{i_1<\dots<i_k \atop j_1<\dots<j_l}} \left.
-m\left(\left[\bigcap_{p\in\left\{i_1,\dots,i_k\right\}}U_p\right]\bigtriangleup
\left[\bigcap_{p\in\left\{j_1,\dots,j_l\right\}}U_p\right]\right)^{2H}
\right\}
\end{split}
\end{equation}
Let us consider the two following terms in expression (\ref{varinterm}):
\begin{itemize}
\item term in $m(U)^{2H}$
\[
m(U)^{2H}+\sum_k (-1)^k\sum_{i_1<\dots<i_k}m(U)^{2H}
=\sum_{k=0}^n C_n^k\; m(U)^{2H}=0
\]
\item term in $m\left(\bigcap_{p\in\left\{i_1,\dots,i_k\right\}}U_p\right)^{2H}$
\begin{align*}
\sum_k (-1)^k\sum_{i_1<\dots<i_k}
m\left(\bigcap_{p\in\left\{i_1,\dots,i_k\right\}}U_p\right)^{2H}
+\underbrace{
\sum_{k,l} (-1)^{k+l}\sum_{i_1<\dots<i_k \atop j_1<\dots<j_l}
m\left(\bigcap_{p\in\left\{i_1,\dots,i_k\right\}}U_p\right)^{2H}
}_{\displaystyle \sum_k (-1)^k\sum_{i_1<\dots<i_k}
m\left(\bigcap_{p\in\left\{i_1,\dots,i_k\right\}}U_p\right)^{2H}
\sum_{l=1}^n (-1)^l C_n^l}
\end{align*}
therefore this term is equal to
\begin{equation*}
\sum_k (-1)^k\sum_{i_1<\dots<i_k}
m\left(\bigcap_{p\in\left\{i_1,\dots,i_k\right\}}U_p\right)^{2H}
\sum_{l=0}^n (-1)^l C_n^l
=0
\end{equation*} 
\end{itemize}
The two other terms of expression (\ref{varinterm}) give the result.
\fin

The main idea to define a set-indexed extension of the fractional Brownian
motion, was to extend
\[
\forall s,t\in\mathbf{R}_{+};\quad 
E\left[\left( B^H_t-B^H_s \right)^2\right]=|t-s|^{2H}
\]
in
\[
\forall U,V\in\mathcal{A};\quad 
E\left[\left(X_U-X_V\right)^2\right]=m(U\bigtriangleup V)^{2H}
\]
However, it should be more interesting to get
\begin{equation}\label{eqCmeas}
\forall C\in\mathcal{C};\quad E\left[\left(\Delta X_C\right)^2\right]=m(C)^{2H}
\end{equation}
According to lemma \ref{statRNlem}, the set-indexed fractional Brownian motion
$\sifbm^H$ satisfies
\[
\forall U,V\in\mathcal{A};\quad 
E\left[\left(\Delta\sifbm^H_{U\setminus V}\right)^2\right]
=m(U\setminus V)^{2H}
\]
but the property (\ref{eqCmeas}) does not hold.\\
Moreover, we will see that there is no set-indexed process satisfying
(\ref{eqCmeas}) for $H\ne\frac{1}{2}$ (theorem \ref{thCstat}).

\begin{proposition}\label{fBsmeas}
Let $\sheet^H=\left\{ \sheet^H_t;\;t\in\mathbf{R}^N_{+} \right\}$ 
be a fractional Brownian sheet of constant parameter $H$ in every axis.\\
For all $a,b\in\mathbf{R}^N_{+}$ such that $a\prec b$ 
(i.e. $\forall i=1,\dots,N;\; a_i<b_i$), we have
\[
E\left[\left( \Delta\left(\sheet^H\right)_{[a,b]} \right)^2\right]
= m([a,b])^{2H}
\]
and consequently,
For all $a,b,a',b'\in\mathbf{R}^N_{+}$, such that $m([a,b])=m([a',b'])$,
we have
\[
\Delta\left(\sheet^H\right)_{[a,b]}\stackrel{(d)}{=}
\Delta\left(\sheet^H\right)_{[a',b']}
\]
\end{proposition}

\proof
For all $a,b\in\mathbf{R}^N_{+}$ with $a\prec b$, 
as $X$ has stationary increments against translations, we have
\begin{align*}
E\left[\left( \Delta\left(\sheet^H\right)_{[a,b]} \right)^2\right]&=
E\left[\left( \Delta\left(\sheet^H\right)_{[0,b-a]} \right)^2\right]\\
&=E\left[\left( \sheet^H_{b-a} \right)^2\right]\\
&=\prod_{i=1}^N |b_i - a_i|^{2H}\\
&= m([a,b])^{2H}
\end{align*}
As $\Delta\left(\sheet^H\right)$ is a mean-zero Gaussian process, 
the result follows. 
\fin

\begin{remark}
\begin{itemize}
\item If a fractional Brownian sheet of index $H=(H_1,\dots,H_N)$ satisfies
the property of proposition \ref{fBsmeas}, we have $H_1=\dots=H_N$.
\item The L\'evy fractional Brownian motion does not satisfy this property.
\end{itemize}
\end{remark}

\begin{definition}\label{statdef}
A set-indexed process $X$ is said to be $\mathcal{C}_0$-increment stationary
if for all $C, C'\in\mathcal{C}_0$ such that $m(C)=m(C')$, we have
$\Delta X_C\stackrel{(d)}{=} \Delta X_{C'}$.
\end{definition}

In the case of one-parameter fractional Brownian motion, the increment
stationarity property gives an equality between laws of some increment processes.
However, in the set-indexed case, definition \ref{statdef}
does not tell anything about correlation between increments.\\
In section \ref{sectcharact}, we see that a stronger property is not worth
to be considered.

\begin{proposition}\label{propstatio}
The set-indexed fractional Brownian motion $\sifbm^H$ is 
$\mathcal{C}_0$-increment stationary.
\end{proposition}

\proof
For all $C\in\mathcal{C}_0$, there exist $U,V\in\mathcal{A}$ 
where $V\subset U$, such that $C=U\setminus V$.
Then $\Delta\sifbm^H_C=\sifbm^H_U - \sifbm^H_{U\cap V}=\sifbm^H_U - \sifbm^H_V$.
We compute
\begin{align*}
E\left[\left(\Delta\sifbm^H_C\right)^2\right]
&=E\left[\left(\sifbm^H_U-\sifbm^H_V\right)^2\right]\\
&=m(U\bigtriangleup V)^{2H}\\
&=m(C)^{2H}
\end{align*}
Thus, as $\Delta\sifbm^H$ is a Gaussian process, for all $C, C'\in\mathcal{C}_0$ 
such that $m(C)=m(C')$, we have
$\Delta\sifbm^H_C\stackrel{(d)}{=} \Delta\sifbm^H_{C'}$.
\fin

\begin{remark}
In the proof of proposition \ref{propstatio}, we saw that
\[
\forall C\in\mathcal{C}_0;\quad
E\left[\left(\Delta\sifbm^H_C\right)^2\right]=m(C)^{2H}
\]
However, in general 
\[
E\left[\Delta\sifbm^H_C .\Delta\sifbm^H_{C'}\right]\neq
\frac{1}{2}\left[m(C)^{2H}+m(C')^{2H}-m(C\bigtriangleup C')^{2H}\right]
\]
In fact, it can be shown that for $C=U\setminus V$ where $V\subset U$, and
$C'=U'\setminus V'$ where $V'\subset U'$, and $U,V,U',V'\in\mathcal{A}$
\[
E\left[\Delta\sifbm^H_C .\Delta\sifbm^H_{C'}\right]=\frac{1}{2}
\left[m(U\bigtriangleup V')^{2H}+m(V\bigtriangleup U')^{2H}
-m(U\bigtriangleup U')^{2H}-m(V\bigtriangleup V')^{2H}\right]
\]
\end{remark}

\subsection{Self-similarity}
To study a set-indexed version of the notion of self-similarity 
for a set-indexed process, we need some assumptions about the set $\mathcal{A}$.
 
We suppose that $\mathcal{A}$ is provided with the operation of a non trivial 
group $G$ that can be extended satifying 
\begin{align}\label{opG}
\forall U,V\in\mathcal{A}, \forall g\in G;\quad
&g.(U\cup V)=g.U \cup g.V\\
&g.(U\setminus V)=g.U \setminus g.V \nonumber
\end{align}
and assume
there exists a non constant function $\mu : G\rightarrow\mathbf{R}^{*}_{+}$
\begin{equation}\label{mu}
\forall U\in\mathcal{A}, \forall g\in G;\quad
m(g.U)=\mu(g).m(U)
\end{equation}

\begin{remark}
We can see easily that $\mu$ is an group-homomorphism from $G$ into the
multiplicative group $\mathbf{R}_{+}\setminus\left\{0\right\}$.
\end{remark}

\begin{example}
In the case of $\mathcal{T}=\mathbf{R}^N_{+}$ and 
$\mathcal{A}=\left\{ [0,t];\; t\in\mathbf{R}^N_{+} \right\} \cup
\left\{\emptyset\right\}$, we can consider 
the multiplication by elements of $\mathbf{R}_{+}$
\begin{equation*}
\forall g\in\mathbf{R}_{+}, \forall t\in\mathbf{R}^N_{+};\quad
g.[0,t]=[0,g.t]
\end{equation*}
Moreover,
\begin{equation*}
\forall g\in\mathbf{R}_{+}, \forall t\in\mathbf{R}^N_{+};\quad
m(g.[0,t])=g^N m([0,t])
\end{equation*}
\end{example}

The following result will be useful in the next section.

\begin{lemma}
Under the assumptions about the group $G$, the cardinal of $G$ is not finite.
\end{lemma}

\proof
As  the function $\mu$ is not constant, 
there exists $g\in G$ such that $\mu(g)> 1$ 
(take $\tilde{g}$ s.t. $\mu(g)\neq 1$ and then 
$g=\tilde{g}$ or $g=\tilde{g}^{-1}$).\\
For all integer $n$, we have $\mu(g^n)=\left[\mu(g)\right]^n$.
If $G$ is finite, the set $\left\{g^n;\;n\in\textbf{N}\right\}$ is finite, 
which is in conflict with 
$\lim_{n\rightarrow\infty}\mu(g^n)=\infty$.
\fin

\begin{definition}
A set-indexed process $X=\left\{X_U;\; U\in\mathcal{A}\right\}$ is said to be
self-similar of index $H$, if there exists a group $G$ which operates on
$\mathcal{A}$, and satisfies (\ref{opG}) and (\ref{mu}),
such that for all $g\in G$,
\begin{equation}
\left\{X_{g.U};\; U\in\mathcal{A} \right\}
\stackrel{(d)}{=} \left\{\mu(g)^H .X_U;\; U\in\mathcal{A} \right\}
\end{equation}

\end{definition}

\begin{proposition}
Assuming the existence of a group $G$ which operates on $\mathcal{A}$, 
and satisfies (\ref{opG}) and (\ref{mu}),
the set-indexed fractional Brownian motion $\sifbm^H$ is self-similar of index 
$H$.
\end{proposition}

\proof
Let $g$ be an element of the group $G$.
For all $U,V \in\mathcal{A}$, we have
\begin{align*}
E\left[\sifbm^H_{g.U} \sifbm^H_{g.V}\right]=
\frac{1}{2}\left[m(g.U)^{2H}+m(g.V)^{2H}-m(g.U\bigtriangleup g.V)^{2H}\right]
\end{align*}
As $g.(U\bigtriangleup V)=g.U\bigtriangleup g.V$, we get
\begin{align*}
E\left[\sifbm^H_{g.U} \sifbm^H_{g.V}\right]&=
\frac{\mu(g)^{2H}}{2}\left[m(U)^{2H}+m(V)^{2H}-m(U\bigtriangleup
V)^{2H}\right]\\
&=\mu(g)^{2H} E\left[\sifbm^H_{U} \sifbm^H_{V}\right]
\end{align*}
Therefore, the two mean-zero Gaussian processes 
$\left\{\sifbm^H_{g.U};\; U\in\mathcal{A} \right\}$ and\\
$\left\{\mu(g)^H .\sifbm^H_U;\; U\in\mathcal{A} \right\}$ have the same law.
\fin

\section{Pseudo-characterisation of sifBm}\label{sectcharact}
Recall that fractional Brownian motion is the only mean-zero Gaussian
process which is self-similar and has stationary increments.
In the same way, the only multiparameter mean-zero Gaussian process which is
self-similar and whose increments are stationary in the strong sense (under
isometries of $\mathbf{R}^N$), is the L\'evy fractional Brownian motion
(\cite{taqqu}).\\
In the case of set-indexed processes, there is not such a characterisation.
However, the two properties of self-similarity and stationarity of increments
characterise the covariance function of the process between any $U$ and $V$ 
such that $U\subseteq V$.

\begin{proposition}\label{characprop}
Suppose the existence of a non trivial group $G$ operating on $\mathcal{A}$, 
and satisfying (\ref{opG}) and (\ref{mu}).
Moreover, assume the fonction $\mu$ is surjective.

Let $X=\left\{ X_U;\;U\in\mathcal{A} \right\}$ be a set-indexed process
satisfying the two following properties
\begin{enumerate}
\item self-similarity of index $H$ with respect to $G$,
\item $\mathcal{C}_0$-increment stationarity
\end{enumerate}
Then, the covariance function between two subsets $U$ and $V$ such that 
$U\subseteq V$ is
\begin{equation}
E\left[ X_U X_V \right]=K\;
\left[ m(U)^{2H}+m(V)^{2H}-m(V\setminus U)^{2H} \right]
\end{equation}
\end{proposition}

\proof
Let $U_0$ be a non $m$-null fixed element of $\mathcal{A}$.
For all $U,V\in\mathcal{A}$ such that $U\subset V$, as $\mu$ is surjective,
there exists $g\in G$ 
such that $\mu(g)=\frac{m(V\setminus U)}{m(U_0)}$, i. e. 
$m(g.U_0)=m(V\setminus U)$. 
Then using $\mathcal{C}_0$-increment stationarity property, we have
\begin{align*}
E\left[\left(X_V-X_U\right)^2\right]
&=E\left[\left(\Delta X_{V\setminus U}\right)^2\right]\\
&=E\left[\left(\Delta X_{g.U_0}\right)^2\right]
\end{align*}
As $g.U_0\in\mathcal{A}$, we have $\Delta X_{g.U_0}=X_{g.U_0}$ and by
self-similarity
\begin{align}\label{covinc}
E\left[\left(X_V-X_U\right)^2\right]
&=E\left[\left(X_{g.U_0}\right)^2\right]\nonumber\\
&=\left[\mu(g)\right]^{2H} E\left[\left(X_{U_0}\right)^2\right]\nonumber\\
&=\left[m(V\setminus U)\right]^{2H} \frac{E\left[\left(X_{U_0}\right)^2\right]}
{m(U_0)^{2H}}
\end{align}
The result follows from (\ref{covinc}).
\fin

\begin{remark}
If $G=\left(\textbf{Q}^{*},.\right)$,
$\mathcal{A}=\left\{[0,t];\;t\in\mathbf{R}^2_{+}\right\} \cup
\left\{\emptyset\right\}$ and $m$ is the Lebesgue measure, we have
\[
\forall q\in\textbf{Q}^{*}, \forall t\in\mathbf{R}^2_{+};\quad
m(q.[o,t])=q^2.m([o,t])
\]
Then, the function $\mu$ is $q\mapsto q^2$, which is not surjective.
The previous result does not hold in this case.
\end{remark}

\begin{remark}
Proposition \ref{characprop} shows that our set-indexed extension of the
fractional Brownian motion is very natural provided the two properties of
self-similarity and stationarity of the increments.

However, if there exists a mean-zero Gaussian process with covariance function
\begin{equation}\label{2ndDef}
E\left[ X_U X_V \right]=\frac{1}{2}
\left[ m(U)^{2H}+m(V)^{2H}-m(V\setminus U)^{2H}-m(U\setminus V)^{2H} \right]
\end{equation}
it satisfies proposition \ref{characprop} as well.
\end{remark}

To determine completely the covariance function of a self-similar,
$\mathcal{C}_0$-increment stationary, set-indexed process, we need assumptions 
about $E\left[\Delta X_{U\setminus V}.\Delta X_{V\setminus U}\right]$, where
$U,V\in\mathcal{A}$.\\
For all $U,V\in\mathcal{A}$, we have 
$\Delta X_{U\setminus V}=X_U - X_{U\cap V}$ and
$\Delta X_{V\setminus U}=X_V - X_{U\cap V}$. Then,
\begin{align*}
E\left[\Delta X_{U\setminus V}.\Delta X_{V\setminus U}\right]&=
E\left[X_U.X_V\right]-E\left[X_U.X_{U\cap V}\right]
-E\left[X_V.X_{U\cap V}\right]+E\left[X_{U\cap V}\right]^2\\
&=E\left[X_U.X_V\right]-\frac{1}{2}
\left[ m(U)^{2H}+m(V)^{2H}-m(V\setminus U)^{2H}-m(U\setminus V)^{2H} \right]
\end{align*}
Particularly, assuming the independance of $\Delta X_{U\setminus V}$
and $\Delta X_{V\setminus U}$, is equivalent to (\ref{2ndDef}), provided that
such a process $X$ exists.
\vsp

The property of $\mathcal{C}_0$-increment stationarity seems too weak to
characterize completely the covariance of a self-similar process.
It can be tempting to define a process which would have a stronger kind of
increment stationarity. For instance, does it exist a self-similar process which
satisfy $E\left[\Delta X_C\right]^2=m(C)^{2H}$ for all $C\in\mathcal{C}$?
The following important result gives a negative answer to the problem of 
existence of such processes.

\begin{theorem}\label{thCstat}
Suppose the existence of a non trivial group $G$ operating on $\mathcal{A}$, 
and satisfying (\ref{opG}) and (\ref{mu}).
Moreover, assume the fonction $\mu$ is surjective and that there exist at least 
two incomparable sets for the partial order $\subset$.

The only Gaussian set-indexed process $X=\left\{ X_U;\;U\in\mathcal{A} \right\}$
such that
\begin{equation}\label{eqCstat}
\forall C\in\mathcal{C};\quad 
E\left[\left( \Delta X_C \right)^2\right] = K.m(C)^{2H}
\end{equation}
where $K>0$ and $H\in (0,1)$, is the set-indexed Brownian motion.
\end{theorem}

\proof
Let $X$ be a set-indexed Gaussian process satisfying (\ref{eqCstat}).\\
First, we can see that $X$ is $\mathcal{C}_0$-increment stationary. Moreover,
for all $U\in\mathcal{A}$ and $g\in G$, we have
$E\left[\left(X_{g.U}\right)^2\right]=K.\mu(g)^{2H}.m(U)^{2H}=\mu(g)^{2H}.E\left[\left(X_U\right)^2\right]$.
Then, as $X$ is Gaussian, we conclude that $X$ is self-similar.
Proposition \ref{characprop} implies that for all $U,V\in\mathcal{A}$, 
such that $U\subseteq V$, 
\begin{equation}\label{eqcovth}
E\left[X_U.X_V\right]=K\left[m(U)^{2H}+m(V)^{2H}
-m(V\setminus U)^{2H}\right]
\end{equation}

For all $U_1$ and $U_2$ in $\mathcal{A}$ such that $U_1\nsubseteq U_2$ and
$U_2\nsubseteq U_1$, let us consider $U\in\mathcal{A}$ such
that $U_1\subset U$ and $U_2\subset U$.
The subset of $\mathcal{T}$, $C=U\setminus (U_1\cup U_2)$ belongs to
$\mathcal{C}$ and
$\Delta X_C=X_U-X_{U_1}-X_{U_2}+X_{U_1\cap U_2}$.
Then we have
\begin{align*}
E\left[\left(\Delta X_C\right)^2\right]=&E\left[\left(X_U\right)^2\right]
+E\left[\left(X_{U_1}\right)^2\right]
+E\left[\left(X_{U_2}\right)^2\right]+E\left[\left(X_{U_1\cap U_2}\right)^2\right]\\
&-2E\left[X_U.X_{U_1}\right]-2E\left[X_U.X_{U_2}\right]
+2E\left[X_U.X_{U_1\cap U_2}\right]\\
&+2E\left[X_{U_1}.X_{U_2}\right]-2E\left[X_{U_1}.X_{U_1\cap U_2}\right]
-2E\left[X_{U_2}.X_{U_1\cap U_2}\right]
\end{align*}
Hence
\begin{align*}
2E\left[X_{U_1}.X_{U_2}\right]&=E\left[\left(\Delta X_C\right)^2\right]
-E\left[\left(X_U\right)^2\right]-E\left[\left(X_{U_1}\right)^2\right]
-E\left[\left(X_{U_2}\right)^2\right]-E\left[\left(X_{U_1\cap U_2}\right)^2\right]\\
&+2E\left[X_U.X_{U_1}\right]+2E\left[X_U.X_{U_2}\right]
-2E\left[X_U.X_{U_1\cap U_2}\right]\\
&+2E\left[X_{U_1}.X_{U_1\cap U_2}\right]
+2E\left[X_{U_2}.X_{U_1\cap U_2}\right]
\end{align*}
Using (\ref{eqcovth}), we get
\begin{align}\label{eqcovU1U2}
2E\left[X_{U_1}.X_{U_2}\right]=&K\left\{m(U\setminus (U_1\cup U_2))^{2H}
-m(U\setminus U_1)^{2H}-m(U\setminus U_2)^{2H}
+m(U\setminus (U_1\cap U_2))^{2H}\right\}\nonumber\\
&+K\left[
m(U_1)^{2H}+m(U_2)^{2H}-m(U_1\setminus U_2)^{2H}-m(U_2\setminus U_1)^{2H}
\right]
\end{align}
Taking $V\in\mathcal{A}$ such that $U\varsubsetneq V$, we get an expression of
$2.E\left[X_{U_1}.X_{U_2}\right]$ different from (\ref{eqcovU1U2}) if
$H\neq\frac{1}{2}$.
\fin

\begin{corollary}
Suppose the existence of a non trivial group $G$ operating on $\mathcal{A}$, 
and satisfying (\ref{opG}) and (\ref{mu}).
Moreover, assume the fonction $\mu$ is surjective and that there exist at least 
two incomparable sets for the partial order $\subset$.

There exists no set-indexed process which is $H$-self-similar (for
$H\ne\frac{1}{2}$) and whose
increments satisfy one of the following
\begin{enumerate}
\item $\mathcal{C}$-increment stationarity
\[
\forall C,C'\in\mathcal{C};\quad m(C)=m(C') \Rightarrow
\Delta X_C\stackrel{(d)}{=}\Delta X_{C'}
\]
\item for all function $f:\mathcal{C}\rightarrow\mathcal{C}$, such that 
$\forall C\in\mathcal{C};\; m(f(C))=m(C)$
\begin{equation*}
\left\{ \Delta X_{f(C)};\;C\in\mathcal{C} \right\} \overset{(d)}{=}
\left\{ \Delta X_{C};\;C\in\mathcal{C} \right\}
\end{equation*}
\end{enumerate}
\end{corollary}

\proof
First, we can see easily that the second property implies 
$\mathcal{C}$-increment stationarity. Then we only need to consider the first
property.\\
Let $U_0$ be a fixed element of $\mathcal{A}$.
For all $C\in\mathcal{C}$, as $\mu$ is surjective, there exists $g\in G$ 
such that $\mu(g).m(U_0)=m(C)$, i.e. $m(g.U_0)=m(C)$.
Therefore, by $\mathcal{C}$-increment stationarity,
\begin{align*}
E\left[\left(\Delta X_C\right)^2\right]&=E\left[\left(X_{g.U_0}\right)^2\right]\\
&=\mu(g)^{2H} E\left[\left(X_{U_0}\right)^2\right]\\
&=m(C)^{2H} \frac{E\left[\left(X_{U_0}\right)^2\right]}{m(U_0)^{2H}}
\end{align*}
By theorem \ref{thCstat}, the result follows.
\fin

\section{Continuity of the sifBm}
The results about the existence of a continuous version of set-indexed processes
are not as simple as processes indexed by $\mathbf{R}_{+}$.
Even in the simple case of the set-indexed Brownian motion, 
if the collection $\mathcal{A}$ is too rich, there does not exist any version
which is continuous on the whole $\mathcal{A}$ (see \cite{adler}).

To study the continuity of a set-indexed process $X$, we have to consider the
behavior of $|X_U-X_V|$ when $U$ and $V$ are close.
In order to do this, we provide $\mathcal{A}$ with some distance.
In the classical case of Gaussian processes, we used to consider the canonical
distance $d^2(U,V)=E\left[\left(X_U-X_V\right)^2\right]$ for $U,V\in\mathcal{A}$ 
(see \cite{adler},\cite{dudley}).
Let us mention two other distances that are also classical: 
\begin{itemize}
\item The measure $m$ on $\mathcal{T}$ induces the pseudo-metric $d_m$ on
$\mathcal{A}$
\[
\forall U,V\in\mathcal{A};\quad d_m(U,V)=m(U\bigtriangleup V)
\]
\item we recall the definition of the Hausdorff metric
$d_{Haus}$\index{$d_{Haus}$}\index{Hausdorff!metric} on
$\mathcal{K}\setminus\emptyset$, the nonempty compact subsets of $\mathcal{T}$
\[
\forall U,V\in\mathcal{K}\setminus\emptyset;\quad
d_{Haus}(U,V)=\inf \{\epsilon >0: U\subseteq V^{\epsilon} \mbox{ and }
V\subseteq U^{\epsilon}\}
\]

\end{itemize}

The notion of continuity depends of the distance considered.
However, if $\left(\mathcal{A},d_m\right)$ 
(resp.$\left(\mathcal{A},d_{Haus}\right)$) is compact, then $d$-continuity and
$d_m$-continuity (resp. $d_{Haus}$-continuity) are equivalent 
(see \cite{adler}).

For any function $x:\mathcal{A}\rightarrow \mathbf{R}$, define
$$\|x\|_{\cal A}=\sup_{A\in{\cal A}} |x(A)|$$
and let
$$B(\mathcal{A} )=\left\{ x:\mathcal{A}\rightarrow\mathbf{R}; 
\|x\|_{\cal A}<\infty \right\}.$$
Let $C\left(\mathcal{A}\right)$ denote the class of functions in 
$B(\mathcal{A})$ which are $d_{Haus}$-continuous on $\mathcal{A}$.
\index{$C({\cal A})$}

Then, studying $d$-continuity, we consider the balls 
$\mathcal{B}(U,\epsilon)=\left\{V\in\mathcal{A};\;d(U,V)<\epsilon\right\}$,
and the metric entropy $D(\bullet;\mathcal{A},d)$ of $(\mathcal{A},d)$, which
gives the smallest number $D(\epsilon;\mathcal{A},d)$ of balls of radius 
$\epsilon>0$ required to cover $\mathcal{A}$.

If $(\mathcal{A},d)$ is totally bounded, i.e. $D(\epsilon;\mathcal{A},d)$ is
finite for all $\epsilon>0$, then under the assumption
\[
\int_0^1 \sqrt{\ln D(\epsilon;\mathcal{A},d)}.d\epsilon <\infty
\]
Dudley's theorem states that
the process $X$ has a continuous modification (see \cite{adler}, \cite{dudley},
\cite{davar}).

\begin{theorem}
Let $\sifbm^H=\left\{\sifbm^H_U;\;U\in\mathcal{A}\right\}$ 
be a set-indexed fractional Brownian motion.\\
The two following statements are equivalent
\begin{enumerate}
\item $\sifbm^H$ is almost surely continuous on $\mathcal{A}$.
\item the set-indexed Brownian motion $\mathbb{W}$ is almost surely continuous on
$\mathcal{A}$.
\end{enumerate}

\end{theorem}

\proof
As
\begin{equation}
\forall U,V\in\mathcal{A};\quad
E\left[\left(\sifbm^H_U-\sifbm^H_V\right)^2\right]
=\left(E\left[\left(\mathbb{W}_U-\mathbb{W}_V\right)^2\right]\right)^{2H}
\end{equation}
the two canonical pseudo-metrics associated to the set-indexed fractional
Brownian motion $\sifbm^H$ and the set-indexed Brownian motion $\mathbb{W}$, 
are equivalent.
\fin

A simple consequence of this result is that the sifBm has a continuous
modification on rectangles of $\mathbf{R}^N_{+}$.

A deeper study of regularity of the set-indexed fractional Brownian motion is
the object of a forthcoming article. In particular, H\"older continuity will be
investigated.

\section{sifBm on increasing paths}
The notion of flows is the key to reducing the proofs of many
of theorems on characterization and weak convergence to a
one-dimensional problem (see \cite{Ivanoff}).

\subsection{Generality on flows}

In general, $\mathcal{A}(u)$ is not closed under countable intersections, so we
will occasionally require a larger class $\tilde{\mathcal{A}}(u)$
\index{$\tilde{\cal A}(u)$}, 
which is the class of countable intersections of  sets in $\mathcal{A}(u)$: 
i.e. $U\in\tilde{\mathcal{A}}(u)$ 
if there exists a sequence $(U_n)_{n\in\textbf{N}}$ in $\mathcal{A}(u)$ 
such that $\bigcap_n U_n= U$.

\begin{definition}
Let $S=[a,b]\subseteq\mathbf{R}$. An increasing function
$f:S \rightarrow\tilde{\mathcal{A}}(u)$ is called a {\it  flow}.
\begin{itemize}
\item A flow $f$ is right-continuous  if
$$f(s)= \bigcap_{v>s} f(v), \forall s\in [a,b),$$
and  $f(b) =\overline{\bigcup_{u<b} f(u)}$.
\index{flow!right-continuous}
\item A flow $f$ is continuous \index{flow!continuous} if it is
right-continuous and
$$f(s)= \overline{\bigcup_{u<s} f(u)}~\forall s\in (a,b).$$
\item A flow $f$ is simple if there exists a finite sequence $(t_0,...t_n)$ 
with
$a=t_0 \leq...\leq t_n=b$ and flows $f_i:[t_{i-1},t_i]\rightarrow\mathcal{A}$,
$i=1,...,n$ such that  for $s\in [t_{i-1},t_i]$, 
$f(s)=f_i(s)\cup \bigcup_{j=1}^{i-1}f_j(t_j)$.
\end{itemize}
\end{definition}

Any  process $X$ indexed by $\mathcal{A}$ can be projected by a simple flow
$f$ onto a process indexed by a subset of $\mathbf{R}$:

\begin{definition}
Let $X$ be an $\mathcal{A}$-indexed process and $f$ a simple flow on $S=[a,b]$.
Then the $S$-indexed process $X^f$ is defined as follows:
$$X^f(s):=X_{f(s)}, ~\forall s\in S.$$
$X^f$ is called the projection of $X$ along $f$.
\end{definition}

In the case that $X$ can be extended to an additive process on
$\tilde{\mathcal{A}}(u)$, $X$ can
be projected by any arbitrary flow according to the preceding definition.

The following lemma shows the importance of the concept of flows for set-indexed
processes.

Let $S(\mathcal{A})$ denote the class of simple continuous flows defined on 
$[0,1]$.

\begin{lemma}[\cite{cime}]
The finite dimensional distributions of an (additive)
$\mathcal{A}$-indexed process $X$ determine and are determined by the finite 
dimensional distributions of the class $\left\{X^f:f\in S(\mathcal{A})\right\}$.
\end{lemma}

In this section, we study the set-indexed fractional Brownian motion on flows.

\subsection{L\'evy fractional Brownian motion and fractional Brownian sheet on
flows}

As a preliminary, let us study the classical cases of L\'evy fractional Brownian
motion and fractional Brownian sheet. \\
Let us consider $\mathcal{T}=\mathbf{R}^N_{+}$ and
$\mathcal{A}=\left\{[0,t];\; t\in\mathbf{R}^N_{+}\right\} \cup
\left\{\emptyset\right\}$.
We can associate to any flow $f$, an increasing function
$\tilde{f}:[0,1]\rightarrow\mathbf{R}^N_{+}$ such that
\[
\forall t\in [0,1];\quad f(t)=[0,\tilde{f}(t)]
\]
\begin{itemize}
\item If $\levy^H$ is a L\'evy fBm of index $H\in (0,1)$,
\begin{equation*}
\forall s,t\in [0,1];\quad 
E\left[\left(\levy^H_{\tilde{f}(t)}-\levy^H_{\tilde{f}(s)}\right)^2\right]=
\| \tilde{f}(t) - \tilde{f}(s) \|^{2H}
\end{equation*}
Then, if $\tilde{f}(t)=t.\alpha$, where $\alpha\in\mathbf{R}^N_{+}$,
$(\levy^H)^{\tilde{f}}$ is a classical fractional Brownian motion, 
otherwise it is not.

\item If $\sheet^H$ is a fBs of index $H\in (0,1)$,
\begin{equation*}
\forall s,t\in [0,1];\quad 
E\left[\sheet^H_{\tilde{f}(s)}.\sheet^H_{\tilde{f}(t)}\right]=\prod_{i=1}^N
\frac{1}{2}\left[ \tilde{f}_i(s)^{2H}+\tilde{f}_i(t)^{2H}
-|\tilde{f}_i(t)-\tilde{f}_i(s)|^{2H} \right]
\end{equation*}
Then, if the function $\tilde{f}$ is a line parallel to one axis of
$\mathbf{R}^N_{+}$, the process $(\sheet^H)^{\tilde{f}}$ is a 
fractional Brownian motion.
However, if $\tilde{f}(t)=t.\alpha$, where $\alpha\in\mathbf{R}^N_{+}$,
\[
\forall s,t\in [0,1];\quad 
E\left[\sheet^H_{\tilde{f}(s)}.\sheet^H_{\tilde{f}(t)}\right]=
\left[ s^{2H}+t^{2H}-|t-s|^{2H} \right]^N
\prod_{i=1}^N\frac{\alpha_i^{2H}}{2}
\]
which is not a fBm.
\end{itemize}
In the two cases of the classical multiparameter extensions of the fractional
Brownian motion, we saw that the projection of the process along a flow, is not
in general a real-indexed fBm.

\subsection{sifBm on flows is a standard fBm}

Our definition for a set-indexed fractional Brownian motion is also justified by
the following proposition.

\begin{proposition}\label{propflow}
Let $\sifbm^H$ be a set-indexed fractional Brownian motion, and $f$ be a
flow on $[0,1]$.
Then the process $(\sifbm^H)^{f}=\left\{\sifbm^H_{f(t)};\; t\in [0,1]\right\}$ 
is a time-changed fractional Brownian motion.
\end{proposition}

\proof
The process $(\sifbm^H)^{f}$ is clearly a mean zero Gaussian process indexed by
$[0,1]$. Moreover, its covariance function can be computed
\begin{align*}
E\left[\sifbm^H_{f(s)}\sifbm^H_{f(t)}\right]=\frac{1}{2}
\left\{ m[f(s)]^{2H}+m[f(t)]^{2H}
-m[f(s)\bigtriangleup f(t)]^{2H} \right\}
\end{align*}
For all $s\leq t$, we have $f(s)\subseteq f(t)$  and then
\begin{align*}
E\left[\sifbm^H_{f(s)}\sifbm^H_{f(t)}\right]&=\frac{1}{2}
\left\{ m[f(s)]^{2H}+m[f(t)]^{2H}-m[f(t)\setminus f(s)]^{2H} \right\}\\
&=\frac{1}{2}\left\{ m[f(s)]^{2H}+m[f(t)]^{2H}
-\left( m[f(t)]-m[f(s)] \right)^{2H} \right\}
\end{align*}
The function $\theta:[0,1]\rightarrow\mathbf{R}_{+}$ such that for all $t\in [0,1]$,
$\theta(t)=m[f(t)]$ is clearly increasing. Thus it defines a time change and
we have
\begin{align}\label{covflow}
\forall s,t\in [0,1];\quad E\left[\sifbm^H_{f(s)}\sifbm^H_{f(t)}\right]=\frac{1}{2}
\left\{ \theta(s)^{2H}+\theta(t)^{2H} - |\theta(t)-\theta(s)|^{2H} \right\}
\end{align}
Then $\left\{\sifbm^H_{f\circ\theta^{-1}(t)};\; t\in\mathbf{R}_{+} \right\}$ 
is a classical fractional Brownian motion.
\fin

Proposition \ref{propflow} allows to identify the self-similarity index of the
sifBm, as the H\"older exponent of the projection along any flow.

Let us recall the definition of the two classical H\"older exponents of a
stochastic process $X$ at $t_0\in\mathbf{R}_{+}$ :
\begin{itemize}
\item the pointwise H\"older exponent
\begin{equation*}
\alphar_X(t_0)=\sup\left\{ \alpha:\;\limsup_{\rho\rightarrow 0}
\sup_{s,t\in \mathcal{B}(t_0,\rho)} \frac{|X_t-X_s|}{\rho^{\alpha}} < \infty
\right\}
\end{equation*}

\item the local H\"older exponent
\begin{equation*}
\tilde{\alphar}_X(t_0)=\sup\left\{ \alpha:\;\limsup_{\rho\rightarrow 0}
\sup_{s,t\in \mathcal{B}(t_0,\rho)} \frac{|X_t-X_s|}{|t-s|^{\alpha}} < \infty
\right\}
\end{equation*}

\end{itemize}

\begin{corollary}
Let $\sifbm^H$ be a set-indexed fractional Brownian motion with self-similarity index
$H$.
The pointwise and local H\"older exponents of the projection $(\sifbm^H)^f$ 
along any flows $f$ at $t_0\in [0,1]$, satisfy almost surely
\begin{align*}
\alphar_{(\sifbm^H)^f}(t_0)&=\left\{
\begin{array}{ll}
\alpha_{\theta}(t_0).H & \textrm{if }\alpha_{\theta}(t_0)<1\\
H & \textrm{otherwise}
\end{array}
\right. \\
\tilde{\alphar}_{(\sifbm^H)^f}(t_0)&=\left\{
\begin{array}{ll}
\tilde{\alpha}_{\theta}(t_0).H & \textrm{if }\tilde{\alpha}_{\theta}(t_0)<1\\
H & \textrm{otherwise}
\end{array}
\right.
\end{align*}
where $\theta$ is the real function such that
$\theta(t)=m\left[f(t)\right]$ ($\forall t\in [0,1]$),
and $\alpha_{\theta}(t_0)$ (resp. $\tilde{\alpha}_{\theta}(t_0)$) is the 
pointwise (resp. local) H\"older exponent of $\theta$ at $t_0$.
 
\end{corollary}

\proof
Let $f$ be a flow and $(\sifbm^H)^f$ the projection of $\sifbm^H$ along $f$.
By (\ref{covflow}), we have
\begin{align}
\forall s,t\in [0,1];\quad
E\left[\left( (\sifbm^H)^f_t - (\sifbm^H)^f_s \right)^2\right]=
|\theta(t) - \theta(s)|^{2H}
\end{align}
\begin{itemize}
\item If $\theta$ is differentiable on $[0,1]$, 
for all $t_0\in [0,1]$ and $\rho>0$, 
\begin{align*}
\forall s,t\in \mathcal{B}(t_0,\rho);\quad
|\theta(t)-\theta(s)|\sim K.|t-s|
\end{align*}
as $\rho$ tends to $0$.
Then,
\begin{align}\label{eqdiff}
\forall s,t\in \mathcal{B}(t_0,\rho);\quad
E\left[\left( (\sifbm^H)^f_t - (\sifbm^H)^f_s \right)^2\right] \sim
K.|t-s|^{2H}
\end{align}
In \cite{ehjlv}, we see that equation (\ref{eqdiff}) implies 
\[
P\left\{ \forall t_0\in [0,1];\;
\alphar_{(\sifbm^H)^f}(t_0)=\tilde{\alphar}_{(\sifbm^H)^f}(t_0)=H 
\right\}=1
\]

\item If $\theta$ is not differentiable in $t_0\in [0,1]$ 
(i.e. if $\tilde{\alpha}_{\theta}(t_0)<1$),
\begin{equation*}
\forall \alpha < \tilde{\alpha}_{\theta}(t_0);\quad
\limsup_{\rho\rightarrow 0}\sup_{s,t\in \mathcal{B}(t_0,\rho)}
\frac{E\left[\left((\sifbm^H)^f_t-(\sifbm^H)^f_s\right)^2\right]}{|t-s|^{2\alpha.H}}=0
\end{equation*}
and 
\begin{equation*}
\forall \alpha > \tilde{\alpha}_{\theta}(t_0);\quad
\limsup_{\rho\rightarrow 0}\sup_{s,t\in \mathcal{B}(t_0,\rho)}
\frac{E\left[\left((\sifbm^H)^f_t-(\sifbm^H)^f_s\right)^2\right]}{|t-s|^{2\alpha.H}}=+\infty
\end{equation*}

then, $\tilde{\alphar}_{(\sifbm^H)^f}(t_0)=\tilde{\alpha}_{\theta}(t_0).H$ almost surely
(see \cite{ehjlv}).\\
In the same way, we get $\alphar_{(\sifbm^H)^f}(t_0)=\alpha_{\theta}(t_0).H$ 
almost surely.
\end{itemize}
\fin

\begin{remark}
As a set-indexed process on a flow only depends on its covariance between
subsets $U$ and $V$ such that $U\subset V$, 
the result stated in proposition \ref{propflow} still holds for
the set-indexed mean-zero Gaussian process defined by (\ref{2ndDef}).
More generally, it holds for all process which is self-similar and
$\mathcal{C}_0$-increment stationary.
Therefore this result lends additional support that such a process could
be called a set-indexed fractional Brownian motion.
\end{remark}

\section{Concluding remarks}
The preceeding discussions permit the following general definition

\begin{definition}
A set-indexed Gaussian process 
which is self-similar of index $H\in (0,1)$ and 
$\mathcal{C}_0$-increment stationary,
is called general set-indexed fractional Brownian motion.
\end{definition}

\begin{corollary}
Let $X=\left\{ X_U;\; U\in\mathcal{A} \right\}$ 
be a general set-indexed fractional Brownian motion. 
Then, 
\begin{enumerate}
\item For all $H\in (0,1/2]$, such a process exists 
(Definition \ref{defsifBm}),
\item for any flow $f$, $X^f$ is a time-changed real-indexed fractional Brownian
motion,
\item the covariance function between two subsets $U$ and $V$ such that 
$U\subseteq V$ is
\begin{equation*}
E\left[ X_U X_V \right]=K\;
\left[ m(U)^{2H}+m(V)^{2H}-m(V\setminus U)^{2H} \right].
\end{equation*}
\end{enumerate}

\end{corollary}

\section*{Acknowledgement}
The authors thank particularly Jacques L\'evy V\'ehel for fruitful discussions about sifBm, and especially for pointing out the counter example of a definition of the sifBm when $H>1/2$.

The second author wishes to thank the first author for his kind invitation, and
Marco Dozzi for preliminary discussions during a stay at Nancy.

\bibliographystyle{plain}
\bibliography{style}

\end{document}